\documentclass{article}
\usepackage{amsfonts}
\usepackage{amsmath}

\input{tcilatex}
\begin{document}

\begin{center}
LARGE ARRAYS OF LINEARLY COUPLED JOSEPHSON JUNCTIONS: A SURVEY
\end{center}

D. G. Aronson

School of Mathematics

University of Minnesota

Minneapolis MN 55455

email: arons001@umn.edu

\bigskip

ABSTRACT: In this note I survey the extensive literature on the dynamics of
large series arrays of identical current biased Josephson junctions coupled
through various shared loads. The equations describing the dynamics of these
arrays are invariant under\ permutation of the junctions so that, in
addition to the usual dynamical systems and numerical methods, group
theoretic methods can also be applied. In practice it is desirable to
operate these circuits at a stable in-phase oscillation. The works
summarized here are devoted to the study of where in parameter space the
in-phase oscillations are stable and how the stability is lost. In
particular, they focus on the variety of states produced through bifurcation
as the in-phase oscillation loses stability. These states include, among
others, discrete rotating waves and semirotors.

\bigskip

Hadley, Beasly \& Wisenfeld [10] studied numerically large series arrays of
current biased Josephson junctions coupled through a variety of shared
loads.\ In practice it is desirable to operate these circuits at a stable
in-phase oscillation. The works summarized here are devoted to the study of
where in parameter space the in-phase oscillations are stable and how the
stability is lost. Here we\ focus on the variety of states produced through
bifurcation as the in-phase oscillation loses stability. A common feature of
the models considered by Hadley \textit{et al.} is that the junctions are
all equally coupled to one another through the load so that the equations
describing their dynamics are invariant under permutation of the junctions.
Thus, in addition to the usual dynamical systems and numerical methods,
group theoretic methods can also be applied [4].

Let $\varphi _{k}(k=1,...,N)$ denote the difference in phase of the
quasiclassical superconducting wavefunctions on the two sides of the $k$-th
junction and let $Q$ denote the current flowing through the load. The
evolution of the $\varphi _{k}$ and $Q$ is governed by the system of
equations:%
\begin{equation}
\beta \overset{..}{\varphi }_{k}+\overset{.}{\varphi }_{k}+\sin (\varphi
_{k})+\overset{.}{Q}=I,  \tag{1}
\end{equation}%
where $\beta $ is the dimensionless measure of the intrinsic capacity of the
junctions, $I$ is the bias current and $\mathfrak{F}$ is an
integro-differential operator which depends on the particular load
considered. Note that the system is invariant under permutation of the
junctions. Here we will concentrate on the load equation:%
\begin{equation}
L\overset{..}{Q}+R\overset{.}{Q}+C^{-1}Q=\frac{1}{N}\tsum_{k=1}^{N}\overset{.%
}{\varphi }_{k},  \tag{2}
\end{equation}%
where $L,R$ and $C$ are, respectively, the inductance, resistance and
capacity of the load.

A \textit{running solution} of the system (1) is an $(N+1)$-tuple%
\begin{equation*}
\left\{ \varphi _{1},\varphi _{2},...,\varphi _{N},Q\right\}
\end{equation*}%
for which there is a minimum $T>0$ such that%
\begin{equation*}
\varphi _{k}(t+T)=\varphi _{k}(t)+2\pi
\end{equation*}%
for $k=1,...,N$ and $Q$ or $\overset{.}{Q}$ is a $T$-periodic function if $%
C<\infty $ or $C^{-1}=0$ respectively. In the sequel we will ignore this
distinction and simply say that $Q$ is $T$-periodic. An in-phase or
symmetric running solution $\left\{ \varphi _{1},\varphi _{2},...,\varphi
_{N},Q\right\} $ is characterized by 
\begin{equation*}
\varphi _{1}=\varphi _{2}=...=\varphi _{N}\text{ .}
\end{equation*}%
If we define $\varphi =\varphi _{1}$ then the system (1),(2) reduces to the
system%
\begin{equation}
\beta \overset{..}{\varphi }+\overset{.}{\varphi }+\sin (\varphi )+\overset{.%
}{Q}=I\text{ \ \ \ \ }L\overset{..}{Q}+R\overset{.}{Q}+C^{-1}Q=\varphi . 
\tag{3}
\end{equation}%
[4] is devoted to the special cases of a pure capacitive load%
\begin{equation}
Q=\frac{3}{N}\sum_{k=1}^{N}\overset{..}{\varphi }_{k}  \tag{4$_{cap}$}
\end{equation}%
and a pure resistive load%
\begin{equation}
Q=\frac{1}{N}\sum_{k=1}^{N}\overset{.}{\varphi }_{k}  \tag{4$_{res}$}
\end{equation}%
with an appropriate scaling suggested by [10]. In these two cases we can
eliminate $Q$ from the system (!),(2) to obtain the system%
\begin{equation}
\beta \overset{..}{\varphi }_{k}+\overset{.}{\varphi }_{k}+\sin (\varphi
_{k})+\frac{\mathcal{A}}{N}\tsum_{j=1}^{N}\overset{.}{\varphi }_{j}+\frac{%
\mathcal{B}}{N}\tsum_{j=1}^{N}\sin (\varphi _{j})=\mathcal{C}I,\text{ \ \ \
\ }k=(1,...,N),  \tag{5}
\end{equation}%
where for a capacitive load%
\begin{equation*}
\mathcal{A=B=}\frac{3}{3+\beta }\text{ \ \ \ \ }\mathcal{C=}\frac{\beta }{%
3+\beta }
\end{equation*}%
and%
\begin{equation*}
\mathcal{A=C=}1,\text{ \ \ \ \ }\mathcal{B=}0
\end{equation*}%
for a resistive load. Symmetric in-phase solutions to (5) are characterized
by $\varphi =\varphi _{1}=\varphi _{2}=...=\varphi _{N}$ . In the capacitive
load case $\varphi $ satisfies the pendulum equation%
\begin{equation*}
(3+\beta )\overset{..}{\varphi }+\overset{.}{\varphi }+\sin (\varphi )=I
\end{equation*}%
while for a resistive load%
\begin{equation*}
\beta \overset{..}{\varphi }+2\overset{.}{\varphi }+\sin (\varphi )=I.
\end{equation*}

The symmetric in-phase oscillations in the pure capacitive and resistive
load cases are analyzed in detail in [4]. It is shown that these
oscillations are born in a homoclinic connection and can only lose stability
through a fixed-point or a period-doubling bifurcation. The $S_{N}$ symmetry
of the equations governing the evolution of an array of $N$ Josephson
junctions is inherited by the Poincar\'{e} map associated with the in-phase
oscillations. A large portion of [4] is devoted to the study of generic $%
S_{N}$ symmetric fixed-point and period-doubling bifurcations for maps. We
classify all possible symmetry breaking period-doubling bifurcations which
can arise in generic period-two states and discuss their stability.. We also
show that generic fixed-point symmetry breaking bifurcations can only
produce unstable fixed points. Note that the group theoretic analysis only
tells us which states cannot occur and which states are possible. For any
specific symmetric system such as the junction arrays considered in [4] the
question of which states actually do occur can only be answered by detailed
analysis of the particular system.

At a period-doubling bifurcation the period-two points of the Poincar\`{e}
map correspond to states where the junctions divide into two or three
groups. Inside each group the junctions oscillate in phase. \ \ When the
number of junctions in each group is different, the oscillations associated
with each group are distinct and their periods are approximately twice the
period 0f the original synchronous solution. When the groups are of the same
size both groups follow the same waveform, but there is a half period phase
shift between groups. If there is a third group of oscillators then,
generically, its size will differ from that of the two equal groups and its
period of oscillation will be half that of the other groups. Thus some
junctions will have period doubled oscillations while others do not.

Determination of the stability of each of the possible period-two points
requires detailed calculations. Essentially the only period-two points which
can be asymptotically stable are those that correspond to two groups of
junctions of approximately equal size. Specifically, if the groups consist
of $k$ and $N-k$ junctions, then generically stability is possible when $%
N/3\leq k\leq N/2.$ Generic $S_{N}$-symmetric fixed-point bifurcations
produce only asymptotically unstable points. Thus generically for each $k$
between $1$ and $N/2$ the fixed-point bifurcations produces fixed points in
which the junctions are divided into two groups of $k$ and $N-k$ junctions.
Within each group the junctions oscillate synchronously. There are no other
fixed points. These theoretic results are born out by detailed numerical
simulations carried out using Doedel's [9] path following software AUTO as
reported in [4]. The numerical results are all consistent with the
theoretical predictions, that is, the solutions found and those not found
are consistent with the theory.

The numerical investigations also led to the discovery of interesting points
on the curve in parameter space along which the homoclinic connections which
give rise to the in-phase solutions occur. These points are \textit{%
homoclinic twist points.} At such points, as the parameters are varied along
the curve of homoclinic points through it, the tangent flow along the
homoclinic trajectory begins to twist vectors in transverse directions.
Three homoclinic twist points were found in the resistive load case and a
unique homoclinic twist point was found in the capacitive load case. The
latter point is of particular interest since it is a codimension-two
bifurcation point where curves of fixed-point and period-doubling
bifurcations intersect and it is called a \textit{homoclinic twist
bifurcation point. }Detailed abstract study of homoclinic twist bifurcations
can be found in [2] and [3]. In [3] we analyze bifurcations occurring in the
vicinity of a homoclinic twist bifurcation point for a generic two-parameter
family of $Z_{2}$ equivariant ODEs in four dimensions. Generically in the
two-parameter space there exists a region of stability of the symmetric
periodic solution bounded by a curve of period-doubling bifurcations and a
curve of pitchfork bifurcations. These bifurcation curves terminate at the
twist point. The period-doubling leads to the creation of a single doubly
winding periodic solution while the pitchfork bifurcation produces two
symmetry-related periodic orbits. Moreover there is a curve in the parameter
space, also terminating at the twist point, where there exists a pair of
symmetry-related homoclinic loops. The periodic orbits born in the
period-doubling and pitchfork bifurcations continue in the parameter space
to the line of the two homoclinic loops and terminate there in an
infinite-period bifurcation. There also exists a branch of doubly winding
homoclinic loops. No homoclinic twist bifurcation points have been found in
the resistive load case. However two the twist points are very close
together and a more refined numerical investigation could show that they
actually coincide. The resulting point would then be a homoclinic twist
bifurcation point.

The dynamics of a pair of identical Josephson junctions coupled through a
shared purely capacitive load are governed by the system%
\begin{eqnarray*}
\beta \overset{..}{\varphi }_{1}+\overset{.}{\varphi }_{1}+\sin \varphi _{1}+%
\frac{C}{2(C+\beta )}(\overset{.}{\varphi }_{1}+\sin \varphi _{1}+\overset{.}%
{\varphi }_{2}+\sin \varphi _{2}) &=&\frac{\beta I}{C+\beta } \\
\beta \overset{..}{\varphi }_{2}+\overset{.}{\varphi }_{2}+\sin \varphi _{2}+%
\frac{C}{2(C+\beta )}(\overset{.}{\varphi }_{1}+\sin \varphi _{1}+\overset{.}%
{\varphi }_{2}+\sin \varphi _{2}) &=&\frac{\beta I}{C+\beta }.
\end{eqnarray*}%
Numerical simulations show that this system possesses a variety of different
running and periodic solutions. Continuation studies using AUTO indicate
that many of these solution branches are generated by a codimension-2
connection which occurs at a particular parameter point which is not a
homoclinic twist point. These branches are described in detail in [1]. There
we also study a general two-parameter system whose properties reflect many
of those found in the numerical studies of the Josephson junction system. In
particular, the model system is assumed to possess an appropriate
codimension-two connection and it is proved that its unfolding generates a
large variety of codimension-1 curves. These results, combined with the
particular symmetry and periodicity properties of junction equations account
for all the numerically observed solution branches in [1]. Indeed, the
theoretical analysis predicted the existence of branches which were not
initially observed. Some of these branches were subsequently found, but most
of them are beyond the reach of numerical simulation. The branch of doubly
winding homoclinic loops found by [4] originates from the homoclinic twist
bifurcation point and terminates in the codimension-2 point found by [1].

In addition to the in-phase or symmetric running solutions, other types of
running solutions occur and have important effects on the overall dynamics
of the full system. Of particular interest are the \textit{discrete rotating
wave solutions}. A running solution is said to be a discrete rotating wave
solution if there exists a running function $\psi $ of period $T$ such that 
\begin{equation}
\varphi _{j}=\psi \left( t-\frac{(j-1)T}{N}\right) \text{ \ \ \ \ }%
j=1,2,...,N  \tag{6}
\end{equation}%
that is, such that every junction follows the same waveform $\psi $ but with
a delay $T/N.$ Discrete rotating wave solutions were first observed
numerically by Hadley \textit{et al}. in [10] and subsequently in [4] were
they are called \textit{pony-on-a-merry-go-round} (POM) solutions. If $N=KM$
then POMs can be formed by grouping the junctions into $M$ blocks of $K$
synchronous junctions with a delay of $T/M$ between successive blocks.
Substitution of (6) in (5) yields the delay-differential equation%
\begin{equation}
\beta \overset{..}{\psi }(t)+\overset{.}{\psi }(t)+\sin \psi (t)+\frac{%
\mathcal{A}}{N}\tsum_{k=0}^{N-1}\overset{.}{\psi }\left( t-\frac{kT}{N}%
\right) +\frac{\mathcal{B}}{N}\tsum_{k=0}^{N-1}\sin \left( \psi (t-\frac{kT}{%
N})\right) =\mathcal{C}I.  \tag{7}
\end{equation}%
This delay equation is special in that the delays $kT/N$ are coupled to the
period $T$ of the running solution $\psi .$ To use global techniques it is
necessary to decouple the period $T$ and introduce a delay parameter $\tau $
leading to the equation%
\begin{equation}
\beta \overset{..}{\psi }(t)+\overset{.}{\psi }(t)+\sin \psi (t)+\frac{%
\mathcal{A}}{N}\tsum_{k=0}^{N-1}\overset{.}{\psi }\left( t-\frac{k\tau }{N}%
\right) +\frac{\mathcal{B}}{N}\tsum_{k=0}^{N-1}\sin \left( \psi (t-\frac{kT}{%
N})\right) =\mathcal{C}I  \tag{8}
\end{equation}%
with $0\leq \tau <N.$ Note that for $\tau =0$ the running solution of (8)
corresponds to the in-phase running solution of (5). In [5] a homotopy and
degree argument is used to extend the known solution of (8) for $\tau =0$ to
prove the existence of the desired solution to (7) for $\tau =T.$

Various authors have studied POMs under different names. Mirollo [14],
Strogatz \& Mirollo [17], Nichols \& Wiesenfeld, among others, call them 
\textit{splay phase solutions}. A compromise position has been agreed upon
and these solutions are now universally known as discrete rotating wave
(DRW) solutions and I will use that terminology from here on. In [14]
Mirollo proves the existence of DRW solutions (which he calls splay phase
solutions) for a pure resistive load in the special case $\beta =0$. Instead
of the homotopy argument used to prove existence in [5] his proof is based
on the Lefschitz trace formula.

So far we have considered junction arrays with only pure capacitive or pure
resistive loads. Arrays with general LCR loads are studied in [8], [12],
[11] and [15]. In [8] Aronson \& Huang consider the system 
\begin{equation}
\beta \overset{.}{\varphi }_{k}+\overset{.}{\varphi }_{k}+\sin \varphi _{k}+%
\overset{.}{Q}=I\text{ for }k=1,2,...,N\text{ \ \ \ \ \ \ }L\overset{..}{Q}+R%
\overset{.}{Q}+C^{-1}Q=\frac{1}{N}\tsum_{j=1}^{N}\overset{.}{\varphi }_{j}, 
\tag{9}
\end{equation}%
where L, R and C are, respectively, the inductance, resistance and
capacitance of the load. Note that $C^{-1}=0$ in case $C=\infty .$ We were
interested in running solutions of this system which follow a single
waveform. A running solution $\left\{ \varphi _{1},\varphi _{2},...,\varphi
_{N},Q\right\} $ of period $T$ is said to be a single wave form solution
(SWFS) if there exists a running function $\phi $ of period $T$ and numbers $%
0=\rho _{1}\leq \rho _{2}\leq ...\leq \rho _{N}<1$, not necessarily
distinct, such that $\varphi _{j}(t)=\phi (t-\rho _{j}T)$ for $j=1,2,...,N$.
An in-phase solution is an SWFS with $\phi =\varphi _{1}$ and $\rho _{j}=0$
for all $j$. A DRW solution is a \ SWFS with $\rho _{j}=(j-1)/N$ for all $j.$

In [8] we define a general class of SWFS for $\beta >0$ which includes
in-phase and DRW solutions and prove their existence. Let $\left\{ \varphi
_{1},\varphi _{2},...,\varphi _{N},Q\right\} $ be a SWFS with waveform $\phi 
$, period $T$ and delays $\{\rho _{1},...\rho _{N}\}.$ Since $\phi $ is a $T$%
-periodic running solution we have%
\begin{equation*}
\beta \overset{..}{\phi }(t)+\overset{.}{\phi }(t)+\sin \phi (t)+\overset{.}{%
Q}(t+T)=I.
\end{equation*}%
Moreover, for each $j\in \{1,2,...,N\}$ it follows from the definition of
SWFS that%
\begin{equation*}
\beta \overset{..}{\phi }(t)+\overset{.}{\phi }(t)+\sin \phi (t)+\overset{.}{%
Q}(t+\rho _{j}T)=I.
\end{equation*}%
Thus%
\begin{equation*}
\overset{.}{Q}(t)=\overset{.}{Q}(t+\rho _{2}T)=...=\overset{.}{Q}(t+\rho
_{N}T)=\overset{.}{Q}(t+T)
\end{equation*}%
and $\phi $ satisfies%
\begin{equation*}
\beta \overset{..}{\phi }(t)+\overset{.}{\phi }(t)+\sin \phi (t)+\overset{.}{%
Q}(t)=I.
\end{equation*}%
Let $\eta (t)$ denote the right hand side of the $Q$-equation in (9). Then%
\begin{equation*}
\eta (t)=\frac{1}{N}\tsum_{j=1}^{N}\overset{.}{\phi }(t-\rho _{j}T)
\end{equation*}%
and it follows that%
\begin{equation*}
\eta (t)=\eta (t+\rho _{2}T)=...=\eta (t+\rho _{N}T)=\eta (t+T).
\end{equation*}

Suppose that $M\leq N$ of the delays $\rho _{j}$ are distinct. Label them%
\begin{equation*}
0=r_{1}<r_{2}<...<r_{M}
\end{equation*}%
and let $l_{1},l_{2},...l_{M}$ denote their multiplicities. Then $l_{j}\in 
\mathbf{N}$ and%
\begin{equation*}
\tsum_{j=1}^{M}l_{j}=N.
\end{equation*}%
For SWFS the pair $\{\phi ,Q\}$ is a solution to the system%
\begin{equation}
\beta \overset{..}{\phi }+\overset{.}{\phi }\sin \phi +\overset{.}{Q}=I\text{
\ \ \ \ \ }L\overset{..}{Q}+R\overset{.}{Q}+C^{-1}Q=\eta (t),  \tag{10}
\end{equation}%
where%
\begin{equation}
\eta (t)=\frac{1}{N}\tsum_{j=1}^{M}l_{j}\overset{.}{\phi }(t-r_{j}T). 
\tag{11}
\end{equation}%
Thus $\phi $ is a $T$-periodic running function, and%
\begin{equation}
\overset{.}{Q}\text{ and }\eta \text{ have periods }r_{2}T,...,r_{M}T,\text{
and }T.  \tag{12}
\end{equation}

Assume that $M\geq 2.$ Since $\overset{.}{Q}$ is $T$-periodic, real valued
and smooth it has a convergent Fourier representation%
\begin{equation*}
\overset{.}{Q}=\tsum_{m}q_{m}e^{2\pi im(t/T)},
\end{equation*}%
where $q_{-m}=\overline{q}_{m}.$ On the other hand, $\overset{.}{Q}$ is also 
$rT$-periodic for every $r\in \{r_{2},...r_{M}\}.$ Therefore we must have%
\begin{equation*}
q_{m}(1-e^{2\pi imr})=0\text{ for all }m\in \boldsymbol{Z.}
\end{equation*}%
If $r$ is irrational then $1-e^{2\pi imr}\neq 0$ for any $m\neq 0$ and it
follows that $\overset{.}{Q}(t)\equiv q_{0}.$ Thus we conclude that if $%
\overset{.}{Q}\neq $constant then all of the delays $r_{2},...,r_{M}$ are
rational numbers. In particular:%
\begin{equation}
\text{\textit{If }}Q\neq \text{\textit{constant} \textit{then} }r_{j}=\frac{%
p_{j}}{K}\text{ \textit{for} }j=2,...,M\mathit{,}\text{\textit{where} }K%
\text{ \textit{is the least common multiple of the denominators of the} }%
r_{j}\text{, }0<p_{2},...,p_{M}\leq K-1\text{ \textit{and} }K\leq M.\text{ 
\textit{If} }K=M\text{ \textit{then} }p_{j}=j-1.  \tag{13}
\end{equation}

Note that it is not necessarily the case that any of the delays in the above
representation are in lowest terms. Write $r_{j}=p_{j}^{\prime }/K_{j}$
where $(rp_{j}^{\prime },K_{j})=1.$ Then $1-e^{2\pi imr_{j}}=0$ only if $%
m\equiv 0$ $(\func{mod}K_{j}).$ Since $K$ is the least common multiple of
the $K_{j}$, it follows that $1-e^{2\pi imr_{j}}\neq 0$ for at least one
value of $j\in \{2,...,M\}$ for every $m\neq 0$ $(\func{mod}K).$ Therefore $%
q_{m}=0$ if $m\neq 0$ $(\func{mod}K).$

In this general case we cannot reduce the existence problem to the study of
a single delay equation as was done in the pure capacitive and resistive
load cases. Instead, to prove existence of a SWFS we introduce a homotopy
parameter $\lambda \in \lbrack 0,1]$ in the system (10)%
\begin{equation*}
\beta \overset{..}{\phi }+\overset{.}{\phi }\sin \phi +\lambda \overset{.}{Q}%
=I\text{ \ \ \ \ \ }L\overset{..}{Q}+R\overset{.}{Q}+C^{-1}Q=\eta (t)
\end{equation*}%
where $\eta (t)$ is given by (11) and the delays satisfy (13). For each $%
\lambda \in \lbrack 0,1]$ we seek a pair $\{\phi ,Q\}$ which satisfies (10),
where $\phi $ is a running solution of period $T(\lambda ),$ and $Q$ and $%
\eta $ satisfy (12) with $T=T(\lambda ).$

For $\lambda =0$ the $\phi $-equation is uncoupled from the $Q$- equation
and is just the simple damped driven pendulum equation It is well known [?]
that there exists a function $I(\beta ):\mathbf{R}^{+}\rightarrow \lbrack
0,1]$ such that for $I>I(\beta )$ the pendulum equation%
\begin{equation*}
\beta \overset{..}{\phi }+\overset{.}{\phi }+\sin \phi =I
\end{equation*}%
possesses a running solution $\phi _{0}$ with period $T_{0}>0.$ Moreover, $%
\overset{.}{\phi }(t)>0$ for all $t$, and is uniquely determined by the
phase condition $\phi _{0}(0)=0.$ Thus \ to start the homotopy we need to
have a solution $Q_{0}$ of the equation%
\begin{equation*}
L\overset{..}{Q}+R\overset{.}{Q}+C^{-1}Q=\frac{1}{N}\tsum_{j=1}^{M}l_{j}\phi
_{0}(t-r_{j}T)
\end{equation*}%
such that $\overset{.}{Q}_{0}$has periods $%
r_{2}T_{0,}r_{3}T_{0,}...,r_{M}T_{0}$ and $T_{0}.$ In the Appendix A of [AH]
it is shown that generically this can only occur if $l_{j}=l$ for all $j\in
\{1,2,..,M\}$, $lM=N$ and 
\begin{equation*}
r_{j}=\frac{j-2}{M}\text{ for }j=1,2,...,M.
\end{equation*}%
Hence we restrict attention to delays which satisfy this necessary condition.

For any positive integers $l,M$ and $N$ such that $lM=N$ we define an $%
(l,M,N)$-SWFS\ of period $T$ \ to be the pair $\{\phi ,Q\}$ such that $\phi $
is a running function of period $T$ and $\overset{.}{Q}$ is a $T/M$-periodic
function such that $\{\phi _{1},\phi _{2},...,\phi _{N},Q\}$is a running
solution to (1),(2) with%
\begin{equation*}
\phi _{j}(t)=\phi \left( t-\frac{h}{M}T\right) \text{ for }j=hl+1,...,(h+1)l%
\text{ and }h=0,1,...,M-1).
\end{equation*}%
Note that an $(N,1,N)$-SWFS is just an in-phase solution, an $(1,N,N)$-SWFS
is an ordinary discrete rotating wave solution and an $(l,M,N)$-SWFS for $%
1<l<M-1$ is a discrete rotating wave solution with $M$ clusters each
consisting of $l$ junctions. The waveform $\phi $ and the load current $%
\overset{.}{Q}$ for an $(l,M,N)$-SWFS satisfy the system%
\begin{equation*}
\beta \overset{..}{\phi }+\overset{.}{\phi }\sin \phi +\overset{.}{Q}=I\text{
\ \ \ \ \ }L\overset{..}{Q}+R\overset{.}{Q}+C^{-1}Q=\frac{l}{N}%
\tsum_{j=1}^{M}\overset{.}{\phi }(t-r_{j}T).
\end{equation*}

In [AH] we prove:

\textbf{Theorem 1:} \textit{Assume that }$M$\textit{\ is a divisor of }$N,$%
\textit{\ }$I\in (1,\infty ),$\textit{\ }$\beta \in (0,\infty ),$\textit{\ }$%
R\in \lbrack 0,\infty ),$\textit{\ }$L\in \lbrack 0,\infty )$\textit{\ and }$%
C\in (0,\infty ]$\textit{\ with }$R+C^{-1}>0.$\textit{\ Then if }$l=N/M$%
\textit{\ there exists an }$(l,M,N)$\textit{-SWFS to (1),(2).}

The theorem is proved by homotopy and degree theory arguments when $%
LC^{-1}=0 $ or $LRC^{-1}>0.$ Limit arguments are used to settle the
remaining cases. Observe that Theorem 1 only asserts the existence of $%
(l,M,N)$-SWFSs and, in particular, does not say anything about the existence
or nonexistence of SWFSs of any other description. Indeed, any SWFS outside
the $(l,M,N)$-class will be invisible to the homotopy argument. However, it
is conjectured that the only SWFSs are those in $(l,M,N)$-class.

In practice, one is usually interested in large arrays of Josephson
junctions and so it is of interest to know how the DRW solutions depend on
the number $N$ of junctions. The existence of a limit of the DRW solutions
and their periods as $N\rightarrow \infty $ is proved in [7] in the pure
capacitive and resistive load cases. The limiting waveform is a solution of
a pendulum equation. In both cases the DRWs are analytic functions of $t$
uniformly in $N$. Thus the rate of convergence is exponential and the
continuum limit provides an excellent approximation even for moderate values
of $N$. Running solutions to (5) are not unique, since any time-shift of a
given solution is again a solution. There is however, an additional source
of non-uniqueness. If a single DRW exists then there are in fact $N!$ of
them obtained by permutation of the indices. Of these, $(N-1)!$ are
genuinely different and the remaining ones can be generated by cyclic
permutations which are equivalent to time-shifts$.$ For general $N$ using
the uniqueness (up to time-shifts) of the running solutions of the pendulum
equation, it is shown in [7[ that the DRWs of (5) are unique up to
permutation and time-shift for sufficiently large $N$.

In [12] Huang\ \& Aronson generalized the results of [7[ on the limiting
behavior as $N\rightarrow \infty $ of DRW solutions to arrays with general $%
LRC$ loads governed by the system (9). Index the waveform $\psi $ and the
charge current $\overset{.}{Q}$ with $N$ and normalize by the condition $%
\psi _{N}(0)=0$ for all $N.$ We prove:

\textbf{Theorem 2:}\textit{\ Suppose }$I>1,\beta >0,C^{-1}\neq 0$\textit{\
or }$R\neq 0$\textit{\ and let }$\{\psi _{N},Q_{N}\}$\textit{\ be a
normalized running solution of (9) with period }$T_{N}>$

\textit{(i) There exist functions }$\psi $\textit{\ and }$Q,$\textit{\ and a
positive constant }$T$\textit{\ such that}%
\begin{equation*}
\lim_{N\rightarrow \infty }T_{N}=T,\text{ }\lim_{N\rightarrow \infty
}Q_{N}(t)=Q(t)\text{ and }\lim_{N\rightarrow \infty }\psi _{N}(t)=\psi (t).
\end{equation*}%
\textit{The convergence is exponential. Moreover all the derivatives of }$%
\psi _{N}$\textit{\ converge to the corresponding derivatives of }$\psi $%
\textit{.}

\textit{(ii) The limiting waveform }$\psi $\textit{\ is a running solution
of the pendulum equation}%
\begin{equation*}
\beta \overset{..}{\psi }+\overset{.}{\psi }+\sin \psi +\frac{2\pi p}{T}=I,%
\text{ where }p=\left\{ 
\begin{tabular}{l}
$0$ if $C^{-1}\neq 0$ \\ 
$R^{-1}$ if $C^{-1}=0$ and $R\neq 0$%
\end{tabular}%
.\right.
\end{equation*}%
\textit{The limiting charge }$Q$\textit{\ or the limiting load current }$%
\overset{.}{Q}$\textit{\ will be constant if }$C<\infty $\textit{\ or }$%
C^{-1}=0$\textit{\ respectively.}

\textit{(iii) There exists an integer }$N_{0}$\textit{\ such that for }$%
N>N_{0}$\textit{\ the running solution of (9) is unique.}

Mirollo \& Rosen also consider arrays with general $LRC$ $\ $loads. Rather
than the homotopy and degree theoretic arguments of [8], they base their
proofs on formulating an\ SWFS\ as a fixed point equation in an appropriate
Hilbert space. Thus, in particular, they do not make any assumptions about
the solution of the pendulum equation and, indeed, derive the existence
directly. Specifically, they prove existence for $I>1$ as well as the
uniqueness of the DRWs for sufficiently large $N$. They also prove

\textit{(iv) There exists an }$I_{0}$\textit{\ depending only on }$L,R$%
\textit{\ and }$C$\textit{\ such \ that the Josephson junction system (9)
admits a unique SWFS for all }$I>I_{0}.$

Finally Mirollo \& Rosen show that, in general, the waveform is not unique.
For this purpose they applied Newton's method to their functional equation
in very precise calculations and found various parameter sets for which
there are multiple waveforms. For example, when%
\begin{equation*}
\beta =4,L=100,R=C=0.01,N=2\text{ and }I=0.7034
\end{equation*}%
they found three distinct waveforms as well as a solution to the pendulum
equation.

Another class of solutions found numerically by Aronson \textit{et al. }in
[4] are the so-called \textit{semirotors}. Here the junctions are split into
two blocks, consisting of $k$ and $N-k$ identical junctions. The junctions
in one block each follow the same $T$-periodic running solution $\psi _{1}$
while each junction in the other block follows a genuinely $T$-periodic
function $\psi _{2}.$ A semirotor is said to be \textit{\ simple} if $\psi
_{1}(t)>0.$ In [6] Aronson, Krupa \& Ashwin study semirotors in the case of
pure resistive and pure capacitive loads. We prove that for sufficiently
large $\beta $ the system (5) has two families of simple semirotors. By
exploiting the symmetries of the system the number of degrees of freedom is
reduced to 4. Numerical studies suggest the form of the solution for small $%
\epsilon =1/\beta $ and the implicit function theorem is used to verify the
actual existence of solution of this form. In the resistive case it is shown
that if $I<(N+k)/N$ for fixed $k\in \{1,...,N-1\}$ there exists $\beta (I)$
such that for $\beta >\beta (I)$ the system (5) has two simple semirotor
solutions \{$\psi _{1}^{j},\psi _{2}^{j}\}$ $(j=1,2)$ depending smoothly on $%
\beta $. The same result holds for the capacitive load case if $I<1$. In
both cases, as $\beta \rightarrow \infty $ the periodic component of the
semirotor shrinks to a single point while the running component approaches a
horizontal line. The stability of the semirotors is also considered in [6].
We compute the asymptotic expansions of the Floquet multipliers of the
linearized system and prove:

\textbf{Theorem 3}: \textit{Let }$\Lambda $\textit{\ be defined by}%
\begin{equation*}
\sin \Lambda =\left\{ 
\begin{tabular}{l}
$\frac{NI}{N+k}$ for a resistive load \\ 
$1$ for a capacitive load%
\end{tabular}%
.\right.
\end{equation*}%
\textit{Then for sufficiently small }$\epsilon $\textit{\ the semirotors
corresponding to }$\cos \Lambda >0$\textit{\ are asymptotically stable and
those corresponding to }$\cos \Lambda <0$\textit{\ are unstable. \ \ \ \ \ \
\ \ \ \ \ \ \ \ \ \ \ \ \ \ \ \ \ \ \ \ \ \ \ \ \ \ \ \ \ \ \ \ \ \ \ \ \ \
\ \ \ \ \ \ \ \ \ \ \ \ \ \ \ \ \ \ \ \ \ \ \ \ \ \ \ \ }Of the two families
of semirotors, one consists of stable solutions and the other of solutions
of saddle--node type. It is conjectured that as $I$ is increased the two
families merge in a saddle-node bifurcation. This conjecture has been
verified numerically (via AUTO) for the case of $N=2$ junctions.

Existence and non-existence of semirotors for arrays with general $LRC$
loads is studied by Huang [11]. She also considers the limiting behavior of
semirotors as $\beta \rightarrow \infty $, but she does not take up the
stability question. She proves:

\textbf{Theorem 4}:\textit{\ (i) Assume that }$R^{2}-4LC^{-1}\neq 0.$\textit{%
\ Then for every fixed }$k\in \{1,2,...,N-1\}$%
\begin{equation*}
\text{when }C^{-1}\neq 0\text{ and }I<1\text{ or when }C^{-1}=0\text{ and }I<%
\frac{k}{NR}
\end{equation*}%
\textit{there exists a function }$\beta (I)$\textit{\ such that the system
(9) has two simple semirotor solutions for every }$\beta >\beta (I).$\textit{%
\ These solutions depend smoothly on }$\beta .$

\textit{(ii) As }$\beta \rightarrow \infty $\textit{\ the phase curve of the
running waveform for a simple semirotor approaches a horizontal line while
the phase curve of the periodic waveform shrinks to a single point.} \ \ \ \
\ \ \ \ \ \ \ \ \ \ The existence proof is based on an extension of the
implicit function theorem method of [6] and does not rule out the
possibility of more than two semirotor solutions.

In the pure resistive and pure capacitive load cases Huang proves a
non-existence result for sufficiently large values of $I$. Specifically:

\textbf{Theorem 5:} \textit{There exist constants }$I(R)$\textit{\ and }$%
I(\beta ,C)$\textit{\ such that for }$I>I(R)$\textit{\ there is no semirotor
solution to the system (5) with pure resistive load and no semirotor
solution to system (5) with pure capacitive load if }$I>I(\beta ,C).$ \ \ \
\ \ \ \ \ \ \ \ \ \ \ \ \ \ \ \ \ \ \ \ \ \ \ \ \ \ \ \ \ \ \ \ \ \ \ \ \ \
\ \ \ \ \ \ \ \ \ \ \ \ \ \ \ \ \ \ \ \ \ \ \ \ \ \ \ \ \ \ \ \ \ \ \ \ \ \
\ \ \ \ \ \ \ \ \ \ \ \ \ \ \ \ \ \ \ \ \ \ \ \ \ \ \ \ \ \ \ \ \ \ \ \ \ \
\ \ \ \ \ \ \ \ \ \ \ \ \ \ \ \ \ \ \ \ \ \ \ \ \ \ \ \ \ \ \ \ \ \ \ \ \ \
\ \ \ \ \ \ \ \ \ \ \ \ \ \ \ \ \ \ \ \ \ \ \ \ \ \ \ \ \ \ \ \ \ \ \ \ \ \
\ \ \ \ \ \ \ \ \ \ \ \ \ \ \ \ \ \ \ \ \ \ \ \ \ \ \ \ \ \ \ \ \ \ \ \ \ \
\ \ \ \ \ \ \ \ \ \ \ \ \ \ \ \ \ \ \ \ \ \ \ \ \ \ \ \ \ \ \ \ \ \ \ \ \ \
\ The main thrust of the proof is showing that the inequalities guarantee
that all of the $\phi _{j}(t)>0$ for large enough $t.$ Actually%
\begin{equation*}
I(R)=\frac{2}{R}+1\text{ and }I(\beta ,C)=1+\frac{2C}{\beta +C}
\end{equation*}%
suffice. Huang conjectures that a similar result holds for general $LRC$
loads.

Everything which we have discussed so far has concerned arrays of identical
Josephson junctions, however it is practically impossible to produce exactly
identical junctions. In [10] Hadley, Beasly \& Wiesenfeld briefly consider
non-identical arrays subject to external noise and present the results of
small-scale numerical simulations. Their simulations of 100 junctions show
that the in-phase solution remains essentially stable when modest junction
mismatches and thermal noise are included. Dhamala \& Wiesenfeld [9] present
a more detailed analysis and show (numerically) the onset of frequency
locking as a function of the statistical disorder in the bias current.
Wiesenfeld, Colet \& Strogatz [20] obtain very precise results on phase
locking for zero-capacity junctions by mapping the junction equations onto
the equations of Sahagachi \& Kuramoto [18] in the limit of weak disorder
and weak coupling. The latter equations are explicitly solvable. \ \ \ \ \ \
\ \ \ \ \ \ \ \ \ \ \ \ \ \ \ \ \ \ \ \ \ \ \ \ \ \ \ \ \ \ \ \ \ \ \ \ \ \
\ \ \ \ \ \ \ \ \ \ \ \ \ \ \ \ \ \ \ \ \ \ \ \ \ \ \ \ \ \ \ \ \ \ \ \ \ \
\ \ \ \ \ \ \ \ \ \ \ \ \ \ \ \ \ \ \ \ \ \ \ \ \ \ \ \ \ \ \ \ \ \ \ \ \ \
\ \ \ \ \ \ \ \ \ \ \ \ \ \ \ \ \ \ \ \ \ \ \ \ \ \ \ \ \ \ \ \ \ \ \ \ \ \
\ \ \ \ \ \ \ \ \ \ \ \ \ \ \ \ \ \ \ \ \ \ \ \ \ \ \ \ \ \ \ \ \ \ \ \ \ \
\ \ \ \ \ \ \ \ \ \ \ \ \ \ \ \ \ \ \ \ \ \ \ \ \ \ \ \ \ \ \ \ \ \ \ \ \ \
\ \ \ \ \ \ \ \ \ \ \ \ \ \ \ \ \ \ \ \ \ \ \ \ \ \ \ \ \ \ \ \ \ \ \ \ \ \
\ \ \ \ \ \ \ \ \ \ \ \ \ \ \ \ \ \ \ \ \ \ \ \ \ \ \ \ \ \ \ \ \ \ \ \ \ \
\ \ \ \ \ \ \ \ \ \ \ \ \ \ \ \ \ \ \ \ \ \ \ \ \ \ \ \ \ \ \ \ \ \ \ \ \ \
\ \ \ \ \ \ \ \ \ \ \ \ \ \ \ \ \ \ \ \ \ \ \ \ \ \ \ \ \ \ \ \ \ \ \ \ \ \
\ \ \ \ \ \ \ \ \ \ \ \ \ \ \ \ \ \ \ \ \ \ \ \ \ \ \ \ \ \ \ \ \ \ \ \ \ \
\ \ \ \ \ \ \ \ \ \ \ \ \ \ \ \ \ \ \ \ \ \ \ \ \ \ \ \ \ \ \ \ \ \ \ \ \ \
\ \ \ \ \ \ \ \ \ \ \ \ \ \ \ \ \ \ \ \ \ \ \ \ \ \ \ \ \ \ \ \ \ \ \ \ \ \
\ \ \ \ \ \ \ \ \ \ \ \ \ \ \ \ \ \ \ \ \ \ \ \ \ \ \ \ \ \ \ \ \ \ \ \ \ \
\ \ \ \ \ \ \ \ \ \ \ \ \ \ \ \ \ \ \ \ \ \ \ \ \ \ \ \ \ \ \ \ \ \ \ \ \ \
\ \ \ \ \ \ \ \ \ \ \ \ \ \ \ \ \ \ \ \ \ \ \ \ \ \ \ \ \ \ \ \ \ \ \ \ \ \
\ \ \ \ \ \ \ \ \ \ \ \ \ \ \ \ \ \ \ \ \ \ \ \ \ \ \ \ \ \ \ 

\ \ \ \ \ \ \ \ \ \ \ \ \ \ \ \ \ \ \ \ \ \ \ \ \ \ \ \ \ \ \ \ \ \ \ \ \ \
\ \ \ \ \ \ \ \ \ \ \ \ \ \ \ \ \ \ \ \ \ \ \ \ \ \ \ \ \ \ \ \ \ \ \ \ \ \
\ \ \ \ \ \ \ \ \ \ \ \ \ \ \ \ \ \ \ \ \ \ \ \ \ \ \ \ \ \ \ \ \ \ \ \ \ \
\ \ \ \ \ \ \ \ \ \ \ \ \ \ \ \ \ \ \ \ \ \ \ \ \ \ \ \ \ \ \ \ \ \ \ \ \ \
\ \ \ \ \ \ \ \ \ \ \ \ \ \ \ \ \ \ \ \ \ \ \ \ \ \ \ \ \ \ \ \ \ \ \ \ \ \
\ \ \ \ \ \ \ \ \ \ \ \ \ \ \ \ \ \ \ \ \ \ \ \ \ \ \ \ \ \ \ \ \ \ \ \ \ \
\ \ \ \ \ \ \ \ \ \ \ \ \ \ \ \ \ \ \ \ \ \ \ \ \ \ \ \ \ \ \ \ \ \ \ \ \ \
\ \ \ \ \ \ \ \ \ \ \ \ \ \ \ \ \ \ \ \ \ \ \ \ \ \ \ \ \ \ \ \ \ \ \ \ \ \
\ \ \ \ \ \ \ \ \ \ \ \ \ \ \ \ \ \ \ \ \ \ \ \ \ \ \ \ \ \ \ \ \ \ \ \ \ \
\ \ \ \textbf{REFERENCES}

1. \ \ D.G. Aronson, E.J. Doedel \& D.H.Terman.\textit{\ A codimension-two
point associated with coupled Josephson junctions.} Nonlinearity \textbf{10}%
(1997) 1231-1255.

2. \ \ D.G. Aronson, S.A. van Gils \& M. Krupa. \textit{The homoclinic twist
point}. International Series of Numerical Mathematics \textbf{104}(1992)
11-22.

3. \ \ D.G. Aronson, S.A. van Gils \& M. Krupa. \textit{Homoclinic Twist
Bifurcations with }$Z_{2}$\textit{\ Symmetry.} J. Nonlinear Science \textbf{4%
}(1994) 195-219.

4. \ \ D.G. Aronson, Martin Golobitsky \& Martin Krupa.\textit{\ Coupled
arrays of Josephson junctions and bifurcation of maps with }$S_{N}$\textit{\
symmetry}. Nonlinearity \textbf{4}(1991) 861-902.

5. \ \ D.G. Aronson, Martin Golobitsky \& John Mallet-Paret.\textit{\ Ponies
on a merry-go-round in large arrays of Josephson junctions.} Nonlinearity 
\textbf{4}(1991) 903-910.

6. \ \ D.G. Aronson, M. Krupa \& P. Ashwin.\textit{\ Semirotors in the
Josephson junction equations}. J. Nonlinear Science \textbf{6}(1996) 85-103.

7. \ \ D.G. Aronson \& Y.S. Huang. \textit{Limits and Uniqueness of discrete
traveling waves in large arrays of Josephson junctions}. Nonlinearity 
\textbf{7}(1994) 777-804.

8. \ \ D.G. Aronson \& Y.S. Huang. \textit{Single waveform solutions for
linear arrays of Josephson junctions.} Physica D \textbf{101}(1997) 157-177.

9. \ \ M. Dhamala \& K. Wiesenfeld. \textit{Generalized stability law for
Josephson series arrays. }Phys. Letters A \textbf{292}(2002) 269-274.

10. E.J.Doedel.\textit{\ AUTO: a program for bifurcation analysis of
autonomous systems.} Congr. Numer. \textbf{30}(1981) 265-384.

11. P. Hadley, M.R. Beasley \& K. Wiesenfeld. Phase locking of Josephson
junctions series arrays. Phys. Rev. B \textbf{38}(1988) 8712-8719.

12. Ying Sue Huang. \textit{Periodic Solutions to Systems of Globally
Coupled Oscillators}. Fields Institute Communications \textbf{21}(1999),
279-291.

13. Y.S.\ Huang \& D.G. Aronson. Discrete Rotating Wave Solutions for
Systems of Globally Coupled Josephson Junctions. International J. Bif. Chaos 
\textbf{6}(1996) 1789-1797.

14. M. Levi, F.C. Hoppensteadt \& M.L. Miranka.\textit{\ Dynamics of the
Josephson Junction}. Q. Appl. Math. \textbf{36(}1978) 167-198.

15. R.E. Mirollo. \textit{Splay-phase orbits for equivariant flow on tori}.
SIAM J. Math Anal. \textbf{25}(1994) 1176-1180.

16. Rennie Mirollo \& Ned Rosen. \textit{Existence, Uniqueness, and
Nonuniqueness of Single Wave-form Solutions to Josephson Junction Systems.}
SIAM J. Appl. Math. \textbf{60}(2000) 1471-1501.

17. Steve Nichols \& Kurt Wiesenfeld.\textit{\ Ubiquitous neutral stability
of splay-states}. Phys Rev. A \textbf{45}(1992) 8430-8435.

18. H. Sahagachi \& Y. Kuramoto.\textit{\ A Soluble Active Rotator Model
Showing Phase Transition via Mutual Entrainment}. Prog. Theor. Phys. \textbf{%
76}(1986) 576-581.

19. S.H. Strogatz \& R.E. Mirollo. \textit{Splay states in globally coupled
Josephson junction arrays: Analytic prediction of Floquet multipliers}.
Phys. Rev. E \textbf{47}(1993) 220-227.

20. Kurt Wiesenfeld, Pere Colet \& Steven H. Strogatz. \textit{Frequency
locking in Josephson arrays: Connections with the Kuramoto Model.} Phys.
Rev. E \ \textbf{57}(1998) 1563.

\bigskip

\bigskip

\bigskip

\bigskip

\bigskip

\bigskip

\bigskip

\bigskip

\end{document}